\documentclass[twocolumn]{IEEEtran}
\usepackage{mathtools,amsmath,amssymb,amsthm,epsfig,color,empheq,graphicx}
\usepackage{enumerate,algorithm,algorithmic,enumitem}
\usepackage[caption=false,font=footnotesize]{subfig}
\usepackage{wasysym,epstopdf}
\usepackage[table]{xcolor}
\usepackage{arydshln,tikz,balance}
\usetikzlibrary{matrix,decorations.pathreplacing}
\usepackage{dashbox, hyperref}

\DeclareMathOperator{\diag}{dg}

\DeclarePairedDelimiter{\ceil}{\lceil}{\rceil}
\DeclarePairedDelimiter{\floor}{\lfloor}{\rfloor}

\newtheorem{lemma}{Lemma}
\newtheorem{corollary}{Corollary}
\newtheorem{theorem}{Theorem}
\newtheorem{definition}{Definition}
\newtheorem{remark}{Remark}

\newcommand \bv{\mathbf{v}}

\newcommand \bE{\mathbf{E}}

\newcommand \bG{\mathbf{G}}

\newcommand \bI{\mathbf{I}}
\newcommand \bJ{\mathbf{J}}

\newcommand \bY{\mathbf{Y}}

\newcommand \mcG{\mathcal{G}}

\newcommand \mcL{\mathcal{L}}
\newcommand \mcM{\mathcal{M}}
\newcommand \mcN{\mathcal{N}}
\newcommand \mcO{\mathcal{O}}

\newcommand \mcR{\mathcal{R}}

\newcommand \mcT{\mathcal{T}}

\newcommand \bmcO{\bar{\mathcal{O}}}

\newcommand \tmcG{\tilde{\mathcal{G}}}

\newcommand \tbJ{\tilde{\mathbf{J}}}

\newcommand \hbp{\hat{\mathbf{p}}}

\newcommand \hbU{\hat{\mathbf{U}}}

\allowdisplaybreaks

\begin{document}
\title{Smart Inverter Grid Probing for Learning Loads:\\
Part I -- Identifiability Analysis}

\author{
	Siddharth Bhela,~\IEEEmembership{Student Member,~IEEE,}
	Vassilis Kekatos,~\IEEEmembership{Senior Member,~IEEE,} and
	Sriharsha Veeramachaneni
	
	
}


\maketitle
\begin{abstract}
Distribution grids currently lack comprehensive real-time metering. Nevertheless, grid operators require precise knowledge of loads and renewable generation to accomplish any feeder optimization task. At the same time, new grid technologies, such as solar photovoltaics and energy storage units are interfaced via inverters with advanced sensing and actuation capabilities. In this context, this two-part work puts forth the idea of engaging power electronics to probe an electric grid and record its voltage response at actuated and metered buses, to infer non-metered loads. Probing can be accomplished by commanding inverters to momentarily perturb their power injections. Multiple probing actions can be induced within a few tens of seconds. In Part I, load inference via grid probing is formulated as an implicit nonlinear system identification task, which is shown to be topologically observable under certain conditions. The conditions can be readily checked upon solving a max-flow problem on a bipartite graph derived from the feeder topology and the placement of probed and non-metered buses. The analysis holds for single- and multi-phase grids, radial or meshed, and applies to phasor or magnitude-only voltage data. Using probing to learn non-constant-power loads is also analyzed as a special case.
\end{abstract}


\section{Introduction}\label{sec:intro}
Low-voltage distribution grids have been plagued with limited observability, due to limited instrumentation, low investment interest in the past, and their sheer extent~\cite{Baran01}. Traditionally, utility operators monitor distribution grids by collecting measurements infrequently and only from a few critical buses. This mode of operation has been functional due to the under-utilization of distribution grids and the availability of historical data. Nevertheless, with the advent of distributed energy resources (DERs), electric vehicles, and demand-response programs, there is a critical need to reliably estimate the system state and learn non-metered loads to optimally dispatch the grid on a frequent basis (say 20 min). To this end, the communication capabilities of grid sensors together with the actuation and sensing features of power inverters found in solar panels, energy storage units, and electric vehicles could be utilized toward unveiling loads. 

Although estimating loads or the grid state has heavily relied on pseudo-measurements, such measurements may not be available or accurate under the current mode of operation~\cite{KlauberZhu15}, \cite{Gomez17}. On the other hand, the widespread deployment of digital relays, phasor measurement units (PMUs), and inverter-interfaced DERs provide excellent opportunities for improving distribution grid observability \cite{DSSE}, \cite{Gao16}. In addition, regular polling and on-demand reads of customer loads and voltages via smart meters have enhanced the accuracy of distribution system state estimation~\cite{Baran09}, \cite{Gray16}. Ignoring network information, a kernel-based scheme for learning loads is reported in~\cite{Jiafan17}. Since the previous schemes collect data on a hourly basis, they are of limited use for real-time optimization.

Rather than passively collecting grid readings to infer non-metered loads, this works advocates engaging inverters to probe the grid and thus actively collect feeder data. We define \emph{probing} as the technique of perturbing an electric grid for the purpose of finding unknown parameters. The idea of probing has been previously suggested towards estimating the electro-mechanical oscillation modes in power transmission systems \cite{ZhTr08}, \cite{TrPi09}. Perturbing the voltage and/or current of a single inverter has been adopted in the power electronics community to determine the grid-equivalent Thevenin impedance of inverters~\cite{JCBR15}. Moreover, modulating the primary droop control loop of inverters has been recently suggested for learning loads and topologies in direct-current grids~\cite{Scal17}. Graph algorithms and identifiability conditions for recovering feeder topologies using inverter probing data have been devised in~\cite{Cavraro18a}, \cite{Cavraro18b}.

Beyond their standard energy conversion functionality, smart inverters are being utilized for reactive power control and other feeder optimization tasks~\cite{Turitsyn11}. In fact, the grid voltage response to inverter injection changes has been used as a means to solve optimal power flow tasks in a decentralized and/or communication-free fashion; see for example~\cite{FCL13}, \cite{VKZG16}, \cite{Emil17}, \cite{Tang18}, \cite{Arnold18}. Leveraging exactly this voltage response, grid probing attributes smart inverters a third functionality towards monitoring rather than grid control.

The contribution of Part I of this work is on three fronts. First, we formulate our \emph{Probing-to-Learn} (P2L) technique in Section~\ref{sec:problem}. Exploiting the stationarity of non-metered loads during probing and assuming noiseless inverter readings in Part I, the P2L problem is posed as a coupled power flow task. Second, we provide intuitive and easily verifiable graph-theoretic conditions under which probing succeeds in finding non-metered loads under phasor (Section~\ref{sec:phasor}) and non-phasor data (Section~\ref{sec:non-phasor}). Lastly, Section~\ref{sec:SSPF} extends probing to infer non-constant-power (ZIP) loads.

The results of Part I significantly extend our previous work of~\cite{BKV17} in four directions: \emph{i)} The analysis extends non-trivially to multiple rather than only two probing actions; \emph{ii)} Probing setups with voltage magnitude and/or angle data are studied in a unified fashion; \emph{iii)} We use the feeder connectivity to upper bound the number of probing actions beyond which there is no identifiability benefit; and \emph{iv)} A new proving technique generalizes the analysis from radial to meshed grids, thus covering the timely topic of loopy multiphase distribution grids and transmission systems.

Regarding \emph{notation}, column vectors (matrices) are denoted by lower- (upper-) case boldface letters and sets by calligraphic symbols. The cardinality of set $\mathcal{X}$ is denoted by $|\mathcal{X}|$, and its complement by $\bar{\mathcal{X}}$. The operators $(\cdot)^{\top}$ and $(\cdot)^{H}$ stand for (complex) transposition; the floor and ceiling functions are denoted by $\floor{\cdot}$ and $\ceil{\cdot}$; $\diag(\mathbf{x})$ defines a diagonal matrix having $\mathbf{x}$ on its main diagonal; and $\bI_{N}$ is the $N \times N$ identity matrix. The notation $\mathbf{x}_{\mathcal{A}}$ denotes the sub-vector of $\mathbf{x}$ indexed by $\mathcal{A}$; and $\mathbf{X}_{\mathcal{A},\mathcal{B}}$ is the matrix obtained by sampling the rows and columns of $\mathbf{X}$ indexed respectively by $\mathcal{A}$ and $\mathcal{B}$.

\section{Grid Probing}\label{sec:problem}
Albeit not every bus is metered in a distribution grid, some buses are equipped with sensors recording voltage magnitudes and/or angles, actual powers, and power factors. Moreover, the power injections in solar panels and energy storage devices can be instantly controlled using advanced power electronics. Building on the physical law that perturbing power injections at different buses is reflected on voltage changes across the grid, the key idea here is to engage power electronics to probe the grid with the purpose of learning non-metered loads. 

To formally describe grid probing, let us briefly review a feeder model. Consider a feeder represented by a graph $\mathcal{G}=(\mathcal{N}^+,\mathcal{L})$ where the nodes in $\mcN^+:=\{0,\ldots,N\}$ correspond to buses, and the edges in $\mathcal{L}$ to distribution lines. Let $\mathbf{Y}:= \mathbf{G} + j\mathbf{B}$ be the grid bus admittance matrix and $\mathbf{G}$ (resp. $\mathbf{B}$) be the bus conductance (resp. susceptance) matrix. By definition, the entries $B_{nm}$ and $G_{nm}$ for $n\neq m$ are non-zero only if $(n,m)\in\mathcal{L}$. Let us express the voltage phasor at bus $n\in\mathcal{N}^+$ in Cartesian and polar coordinates as
 \[v_n=v_{r,n}+jv_{i,n}=u_n e^{j\theta_n}.\]
The substation is indexed by $n=0$, its voltage remains fixed at $1+j0$, and the remaining buses comprise the set $\mcN$. If $\bv_r:=[v_{r,0}~\cdots~v_{r,N}]^\top$ and $\bv_i:=[v_{i,0}~\cdots~v_{i,N}]^\top$, define the system state as $\bv:=[\bv_r^{\top}~\bv_i^{\top}]^{\top}$. Apparently, for each bus $n\in \mcN^+$, the squared voltage magnitude and the net power injections are quadratic functions of $\bv$, whereas the voltage angle is a trigonometric function of $\bv$~\cite[Ch.~3]{ExpConCanBook}
\begin{subequations} \label{eq:PF}
	\begin{align}
	u_n(\bv)&=u_n^2=v_{r,n}^2 + v_{i,n}^2\label{eq:PFu} \\
	p_n(\mathbf{v})&=v_{r,n} \sum_{m=0}^N \left(v_{r,m}G_{nm}-v_{i,m}B_{nm}\right) \nonumber\\
	&\quad+v_{i,n} \sum_{m=0}^N\left(v_{r,m}B_{nm}+v_{i,m}G_{nm}\right) \label{eq:PFa}\\
	q_n(\mathbf{v})&=v_{i,n} \sum_{m=0}^N \left(v_{r,m}G_{nm}-v_{i,m}B_{nm}\right) \nonumber\\
	&\quad-v_{r,n} \sum_{m=0}^N \left(v_{r,m}B_{nm}+v_{i,m}G_{nm}\right) \label{eq:PFb}\\
	\theta_n(\bv)&= \arctan\left(\frac{v_{r,n}}{v_{i,n}}\right).\label{eq:PFtheta}
	\end{align}
\end{subequations}

With the proliferation of grid sensors and inverters, the distribution grid operator may have access to all four quantities $(u_n,\theta_n,p_n,q_n)$ on a subset of buses. Different from the conventional power flow (PF) setup with PQ and PV buses, we partition $\mathcal{N}^+$ into the subsets:
\begin{itemize}
\item The set $\mathcal{M}$ of metered buses for which $(u_n,\theta_n,p_n, q_n)$ are known and their power injections are possibly controllable. This set includes the substation and buses equipped with smart sensors and/or inverters. Its cardinality is denoted by $M:=|\mcM|$.
\item The set $\mcO$ of non-metered buses where no information is available. Its cardinality is denoted by $O:=|\mcO|$, and apparently, $N+1=M+O$.
\end{itemize}

The inverters interfacing DERs are typically modeled as constant-power generators~\cite{Turitsyn11}, \cite{FCL13}, \cite{Emil17}: Internal control loops can reach setpoints for (re)-active power injections within microseconds. The setpoints should comply with solar irradiance and the rating of the inverter.

\begin{remark}\label{re:pinjection}
We emphasize $p_n+j q_n$ is the net complex injection. If bus $n$ hosts a smart inverter and a non-controllable load, it is henceforth assumed that the operator measures $p_n+j q_n$ and the voltage at the point of common coupling, and controls the complex injection from the inverter. This assumption is reasonable since smart inverters are usually equipped with sensors; e.g., the Pecan Street project measures both the net and inverter injections~\cite{pecandata}.
\end{remark}

Given the feeder topology captured in $\bY$ and the specifications $\{(u_n,\theta_n,p_n,q_n)\}_{n\in\mcM}$, our goal is to recover the power injections at non-metered buses $\{(p_n,q_n)\}_{n\in\mcO}$. Lacking a direct mapping from $\{(u_n,\theta_n,p_n,q_n)\}_{n\in\mcM}$ to $\{(p_n,q_n)\}_{n\in\mcO}$, the problem of finding the non-metered loads boils down to the task of recovering the underlying state $\bv$ first. Collecting the grid data $\{(u_n^t,\theta_n^t,p_n^t,q_n^t)\}_{n\in\mcM}$ at time $t$ and assuming for now these data are noiseless, we get the specifications
\begin{subequations}\label{eq:PF1}
	\begin{align}
	u_n(\bv_t)&=u_n^t, &\forall n\in\mcM\label{eq:PF1u}\\
	\theta_n(\bv_t)&=\theta_n^t, &\forall n\in\mcM\label{eq:PF1t}\\
	p_n(\bv_t)&=p_n^t,~q_n(\bv_t)=q_n^t &\forall n\in\mcM\label{eq:PF1pq}
	\end{align}
\end{subequations} 
which involve $4M$ equations over $2(N+1)$ unknowns. A necessary condition for solving \eqref{eq:PF1} is $4M\geq 2(N+1)$. Since $N+1=M+O$, the condition simplifies to
\begin{equation}\label{eq:cond1}
M\geq O.
\end{equation}
In other words, the metered buses must be at least as many as the non-metered ones. 

\begin{figure*}[h]
\centering
\includegraphics[scale=0.4]{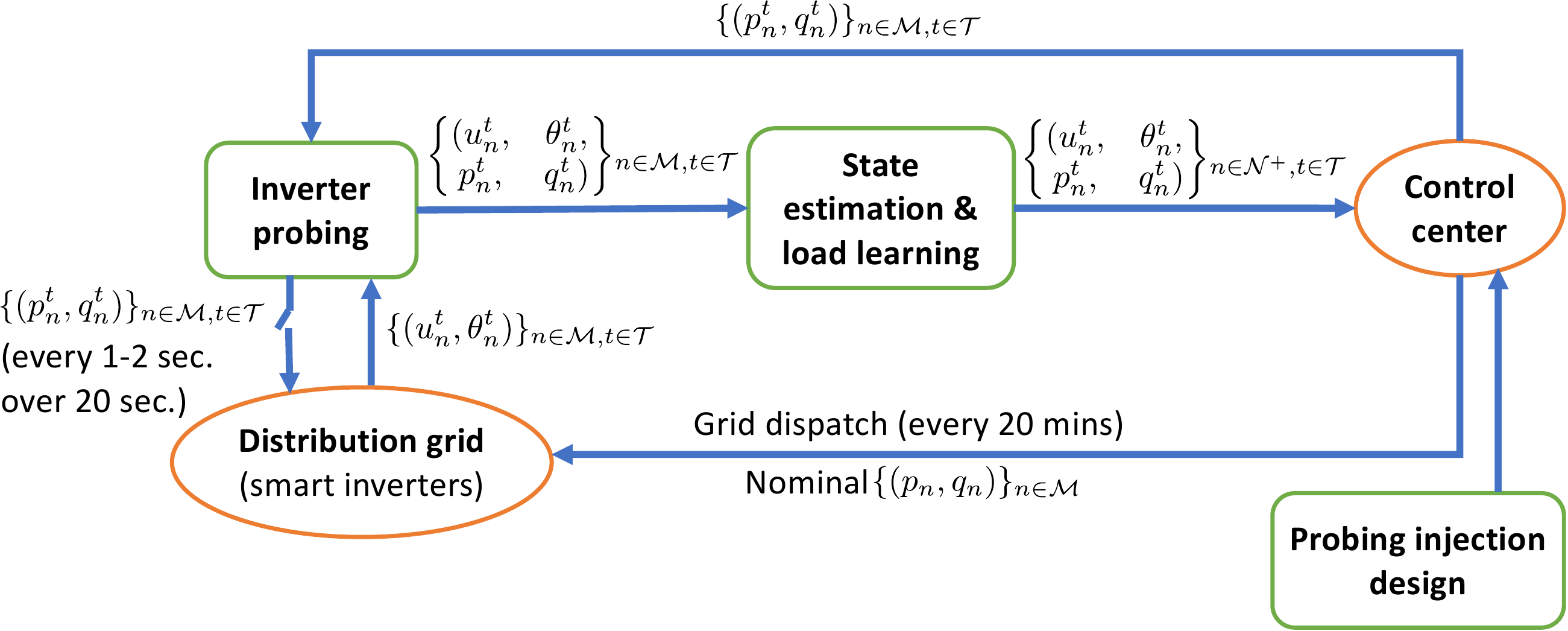}\hspace*{1em}
\includegraphics[scale=0.44]{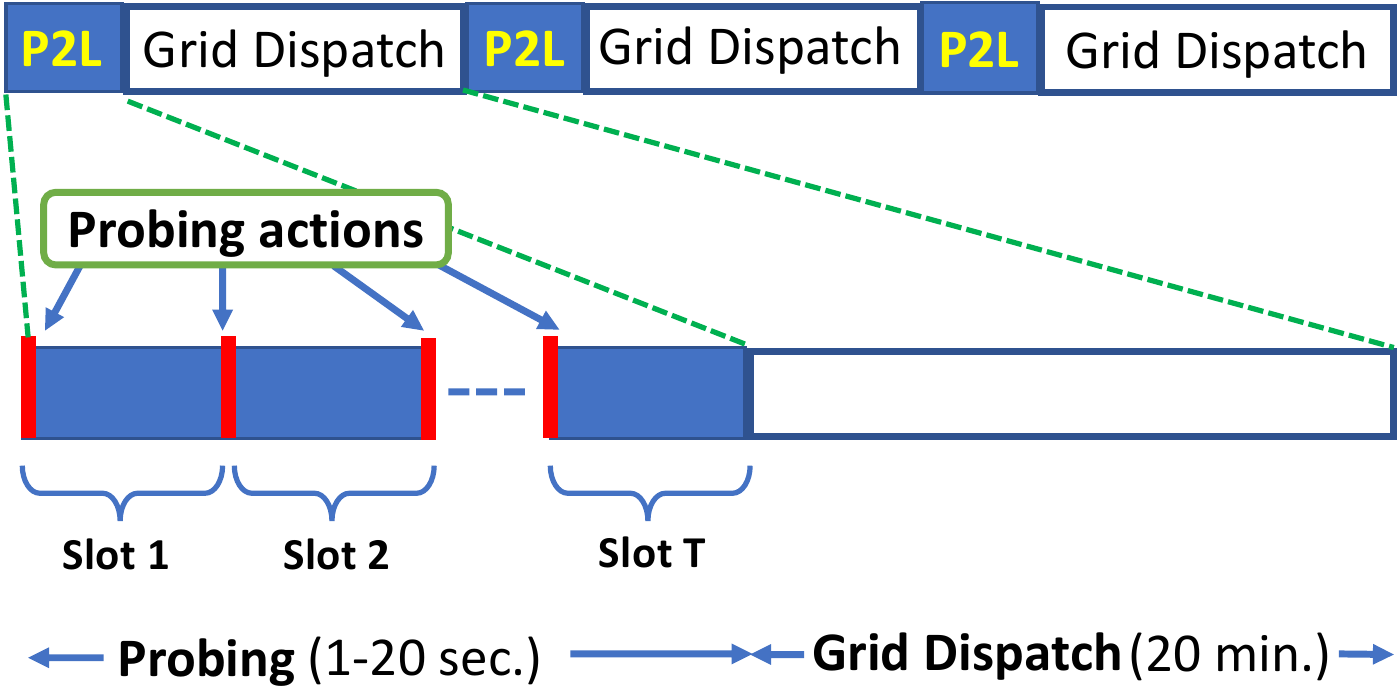}
\caption{Overview of the P2L framework: (a) block diagram depicting the P2L task with phasor data; (b) temporal organization of grid operation.}
\label{fig:flowchart}
\end{figure*}

To relax this condition on $M$, one may consider jointly processing the data $\{(u_n^t, \theta_n^t, p_n^t,q_n^t)\}_{n\in\mcM}$ collected across multiple times $t\in\mcT$ with $\mcT:=\{1,\ldots,T\}$. This approach does not improve the observability of the equations in \eqref{eq:PF1}, simply because the equations are independent over $\mcT$. Moreover, both the $4MT$ equations and the $2(N+1)T$ state variables $\{\bv_t\}_{t=1}^T$ scale with $T$. One way to relate power flow specifications across time is to assume that the non-metered loads remain invariant across $\mcT$, that is
\begin{subequations}\label{eq:PFC}
	\begin{align}
	p_n(\bv_t)&=p_n(\bv_{t+1}), &\forall n\in\mcO,t \in \mcT' \label{eq:PFCp}\\
	q_n(\bv_t)&=q_n(\bv_{t+1}), &\forall n\in\mcO,t \in \mcT' \label{eq:PFCq}
	\end{align}
\end{subequations} 
where $\mcT':=\{1,\ldots,T-1\}$. In this way, we obtain the additional $2O(T-1)$ equations and couple the states $\{\bv_t\}_{t=1}^T$.

Even though there may be an observability advantage in coupling specifications across time, the timespan of $\mcT$ is critical: For non-metered loads to remain unchanged, the timespan of $\mcT$ should be relatively short. But if the duration of $\mcT$ is too short, the metered injections in the buses of $\mcM$ may not change either. In this case, the grid state remains identical over $\mcT$, the scheme degenerates to the setup of \eqref{eq:PF1} for $T=1$, and there is no advantage by coupling specifications. 

At this point, smart inverters come to our rescue: The timespan of $\mcT$ can be made sufficiently short so that the non-metered loads in the buses of $\mcO$ remain invariant over $\mcT$, whereas the power injections from smart inverters vary. The key point here is to couple power flow specifications through what we term \emph{grid probing}. Probing can be accomplished by commanding inverters to change their power injections for one second. An inverter can curtail its solar output; (dis)-charge an energy storage unit; and/or change its power factor. Multiple probing actions can be instructed within tens of seconds. By intentionally perturbing inverter injections, the grid transitions across different states $\{\bv_t\}_{t=1}^T$ depending on the probing injections and non-metered loads. Recording voltages $\{u_n^t,\theta_n^t\}_{t=1}^T$ over $n\in\mcM$ could unveil non-metered loads. 

The metered buses in $\mcM$ can be classified into probing buses and metered but non-controllable buses. Although grid data $(u_n^t, \theta_n^t, p_n^t,q_n^t)$ are collected on both probed and metered buses, the operator can control only the probing buses. To simplify the presentation, we will henceforth assume that all metered buses are probing buses, although the analysis and algorithms apply to the more general setup.

Probing postulates two assumptions on non-metered loads: \emph{a1)} They remain constant throughout $\mcT$; and \emph{a2)} are modeled as of constant power. Assumption \emph{a1)} may be reasonable over the short duration of probing. Regarding \emph{a2)}, one could alternatively adopt a ZIP load model for bus $n\in\mcO$~\cite{Kersting}
\begin{subequations} \label{eq:ZIP0}
	\begin{align}
	-p_n^t(u_n^t)&=\alpha_{p_n} (u_n^t)^2 +\beta_{p_n} u_n^t +\gamma_{p_n} \label{eq:p}\\ 
	-q_n^t(u_n^t)&=\alpha_{q_n} (u_n^t)^2 +\beta_{q_n} u_n^t +\gamma_{q_n} \label{eq:q}.
	\end{align}
\end{subequations}
The parameters $(\alpha_{p_n},\beta_{p_n},\gamma_{p_n})$ correspond to the constant-impedance, constant-current, and constant-power components of active load; likewise $(\alpha_{q_n},\beta_{q_n},\gamma_{q_n})$ for reactive load. A non-metered ZIP load is then described by six rather than two parameters. Moreover, despite the model for load, $n$ does not change across $\mcT$, its power injection $p_n+jq_n$ does change for varying $u_n$. Then, the coupling equations are not valid in the form of \eqref{eq:PFC} anymore. If the ZIP parameters are assumed invariant over $\mcT$, the power flow equations can still be coupled across $\mcT$, yet the identifiability analysis and the associated solvers become perplex. To bypass this complexity, Section~\ref{sec:SSPF} copes with ZIP loads by resorting to single-slot probing. 

Figure~\ref{fig:flowchart} depicts how probing can be incorporated into grid operation: Suppose a utility operates a demand-response program; manages energy storage; or controls smart inverters for reactive power control on a 20-min basis. To solve the optimal power flow problem, the operator needs to know the injections at non-metered buses. To do so, a probing interval lasting few tens of seconds precedes the feeder dispatch. This interval $\mcT$ is divided into $T$ probing slots indexed by $t=1,\ldots,T$. During each probing slot $t$, every inverter $n\in\mcM$ changes its injections to the setpoints $(p_n^t,q_n^t)$ and reads voltage data $(u_n^t)$ or $(u_n^t,\theta_n^t)$. At the end of interval $\mcT$, each inverter $n\in\mcM$ sends the collected data $\{u_n^t\}_{t\in\mcT}$ or $\{u_n^t,\theta_n^t\}_{t\in\mcT}$ back to the utility, and switches its setpoints back to their nominal values. The utility processes the collected data, infers the non-metered loads, and dispatches the grid for the next 20-min period. Some implementation details follow.

\begin{remark}\label{re:coordination}
The probing setpoints $(p_n^t,q_n^t)$ for all $t$ and $n\in\mcM$ are decided by the utility prior to $\mcT$ and communicated to all inverters via two-way communication links. This is to ensure that probing complies with voltage constraints and for improved load estimation accuracy; see Part II. The commanded setpoints are attained by simple PID controllers. Further, the inverters act synchronously along probing slots. Since potential delays may raise synchronization issues, developing protocols where inverters probe asynchronously is of interest.
\end{remark}

\begin{remark}\label{re:topology}
The proposed probing scheme aims at recovering loads assuming the feeder topology is known. The topology includes bus connectivity and line impedances, phase assignments, and the statuses of capacitors and voltage regulators; see Remark~\ref{re:switching}. Although small errors in line impedances and regulator tap settings could be modeled as measurement noise, grid probing is as sensitive to topology errors as power flow equations are. However, probing can be also used for inferring grid topologies and line parameters without knowing non-metered loads~\cite{Cavraro18a}, \cite{Cavraro18b}. Moreover, phase assignments can be inferred from smart meter data; see e.g., \cite{liao2018}. Such techniques could precede P2L to find or calibrate feeder models.
\end{remark}

\begin{remark}\label{re:switching}
Feeders are equipped with voltage-control devices, such as regulators and capacitor banks, which respond to voltage excursions by changing their taps and switching on/off with time delays of around 30-90~seconds~\cite{Kersting}. Since $\mcT$ lasts 20~sec or less, probing is not expected to trigger voltage control actions per se. Nonetheless, there are still chances for these actions to occur during $\mcT$ due to load fluctuations. If the utility does not monitor these devices in real-time or it cannot override their settings during probing, the topology learning techniques of Remark~\ref{re:topology} could be possibly used. Voltage control actions and topology reconfigurations will be ignored in this work. Interestingly though, such actions could be used towards grid probing too.
\end{remark}

Grid probing can be now formally stated as follows.

\begin{definition}[Probing-to-Learn task with phasor data]\label{def:P2L}
Given $\mathbf{Y}$ and probing data $(u_n^t, \theta_n^t, p_n^t,q_n^t)$ for all $n\in\mcM$ and $t\in\mcT$, the probing-to-learn (P2L) task entails solving the equations in \eqref{eq:PF1} for $t\in\mcT$ jointly with the coupling equations in \eqref{eq:PFC}. 
\end{definition}

The P2L task involves $4MT+2O(T-1)$ equations in $2(N+1)T$ unknowns. A necessary condition for solving it is
\begin{equation}\label{eq:cond2}
M\geq \frac{O}{T}
\end{equation}
which coincides with the condition in \eqref{eq:cond1} for $T=1$. For $T\geq 2$ however, it improves upon \eqref{eq:cond1} if probing over multiple time instances is allowed. In \cite{BKV17}, we have derived conditions under which the P2L task recovers non-metered loads for $T=2$. The analysis there was further confined to non-phasor grid data $\{(u_n^t,p_n^t,q_n^t)\}_{n\in\mcM,t \in \mcT}$ and radial grids. The conference work of \cite{BKVCS17} extended the previous claims (without proofs) to meshed networks. Here, we broaden the scope to study the identifiability of the P2L task with phasor data $\{(u_n^t,\theta_n^t, p_n^t,q_n^t)\}_{n\in\mcM}$ over $\mcT$, and show that the analysis with non-phasor data can be seen as a special case of the former.

\section{Identifiability of P2L with Phasor Data} \label{sec:phasor}
As customary in identifiability analysis, data will be assumed noiseless; noisy data are considered in Part II~\cite{BKV18P2}. The relationship between the inputs $\{u_n^t,\theta_n^t, p_n^t,q_n^t\}_{n\in\mcM}$ and the outputs $\{p_n^t,q_n^t\}_{n\in\mcO}$ of the P2L task is implicit since the PF equations involve $\{\bv_t\}_{t=1}^T$ as nuisance variables. Because of this, P2L is tackled in two steps. The first step of finding $\{\bv_t\}_{t=1}^T$ is the challenging one. In the second step, one simply evaluates $(p_n(\bv_t),q_n(\bv_t))$ for all $n\in\mcO$ and $t=1$. For numerical stability, one can recover the unknown injections by averaging as $\frac{1}{T}\sum_{t=1}^T p_n(\bv_t)$ and $\frac{1}{T}\sum_{t=1}^T q_n(\bv_t)$ for all ${n\in\mcO}$. Hence, if the system states $\{\bv_t\}_{t=1}^T$ can be recovered by solving \eqref{eq:PF1} and \eqref{eq:PFC}, the P2L task is deemed successful. 

Granted the P2L equations are non-linear, identifiability can be ensured only within a neighborhood of the nominal $\{\bv_t\}_{t=1}^T$. Upon invoking the inverse function theorem, a necessary and sufficient condition for locally solving P2L is that the Jacobian matrix $\bJ\left(\{\bv_t\}\right)$ related to the nonlinear equations of \eqref{eq:PF1} and \eqref{eq:PFC} is full rank. Because $\bJ\left(\{\bv_t\}\right)$ depends on $\{\bv_t\}$, characterizing its column rank for any $\{\bv_t\}$ is challenging. 

To tackle this issue, we resort to the \emph{generic rank} of a matrix defined as the maximum possible rank attained if the non-zero entries of the matrix are allowed to take arbitrary real values~\cite{Tutte}, \cite{Tucker04}. If the generic rank of an $M \times N$ matrix $\bE$ with $M\geq N$ equals $N$, matrix $\bE$ is said to be of \emph{full generic rank}. The generic rank of a matrix is related to a graph constructed by the sparsity pattern of the matrix, that is the locations of its (non)-zero entries. To explain this link, some graph-theoretic concepts are needed. 

\begin{figure}[t]
\centering
\includegraphics[scale=0.56]{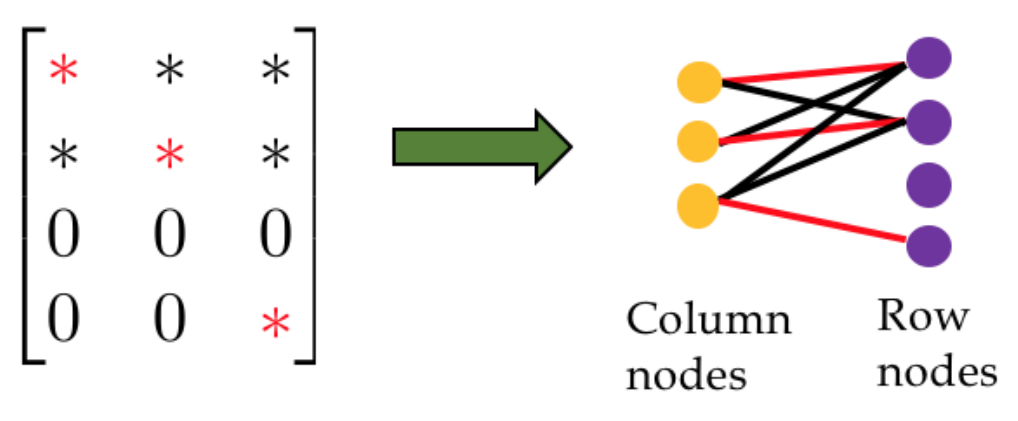}
\caption[font=scriptsize]{The sparsity pattern of $\bE$ and its bipartite graph $\mcG_E$: Column nodes are linked to row nodes depending on the entries of $\bE$. The perfect matching is marked in red. From Lemma~\ref{le:bipartite}, any matrix with this sparsity pattern is generically full rank. Had $E_{4,3}=0$, no perfect matching would exist.}
\label{fig:lemma}
\end{figure} 

A graph $\mcG=(\mcN,\mcL)$ is bipartite if $\mcN$ can be partitioned into disjoint subsets $\mcN_1$ and $\mcN_2$, such that $\mcN=\mcN_1 \cup \mcN_2$, and every $\ell \in \mcL$ connects a node in $\mcN_1$ to a node in $\mcN_2$. A subset of edges $\mcL'\subseteq \mcL$ is termed a perfect matching of $\mcN_1$ to $\mcN_2$, if every vertex in $\mcN_1$ is incident to exactly one edge in $\mcL'$. The degree $\delta_n(\mcG)$ of node $n$ is defined as the number of edges incident to node $n$ in $\mcG$. Given a matrix $\bE\in\mathbb{R}^{M\times N}$, construct a bipartite graph $\mcG_E$ having $M+N$ nodes: Each column of $\bE$ is mapped to a column node and each row of $\bE$ to a row node. An edge runs from the $n$-th column node to the $m$-th row node only if $E_{mn}\neq 0$; see Fig.~\ref{fig:lemma}. Based on $\mcG_E$, we will use next claim.

\begin{lemma}[\cite{Tutte}, \cite{Tucker04}]\label{le:bipartite}
An $M\times N$ matrix $\bE$ has full generic rank if and only if the bipartite graph $\mcG_E$ features a perfect matching from the column nodes to its row nodes.
\end{lemma}

According to Lemma~\ref{le:bipartite} (proved in \cite[Th.~12.10]{Tucker04}), the generic identifiability of P2L relies on the sparsity pattern of $\bJ\left(\{\bv_t\}\right)$. The goal is to match every column node (state) of $\bJ\left(\{\bv_t\}\right)$ to a unique row node (equation). The non-zero entries of $\bJ\left(\{\bv_t\}\right)$ are the available links. 

To characterize the sparsity pattern of $\bJ\left(\{\bv_t\}\right)$, consider the Jacobian matrices $\mathbf{J}^u(\bv)$, $\mathbf{J}^{\theta}(\bv)$, $\mathbf{J}^p(\bv)$, and $\mathbf{J}^q(\bv)$, associated accordingly with the squared voltage magnitudes and voltage angles, and the (re)active power injections over all buses. Matrix $\bJ\left(\{\bv_t\}\right)$ consists of stacked row-sampled submatrices of $\bJ^u(\bv_t)$, $\bJ^{\theta}(\bv_t)$, $\bJ^p(\bv_t)$, and $\bJ^q(\bv_t)$ corresponding to \eqref{eq:PF1} and \eqref{eq:PFC} for $t\in\mcT$. The matrices obtained by selecting the rows of $\bJ^u(\bv_t)$ associated with buses in $\mcM$ and $\mcO$ are respectively denoted by $\bJ^u_\mcM(\bv_t)$ and $\bJ^u_\mcO(\bv_t)$. Similar notation is used for $\bJ^{\theta}(\bv_t)$, $\bJ^p(\bv_t)$, and $\bJ^q(\bv_t)$. Let us define
\begin{subequations}
	\begin{align*}
	\bJ_{\mcM}(\bv_t)&:= \left[\begin{array}{c}
	\bJ^u_\mcM(\bv_t) \\
	\bJ^{\theta}_\mcM(\bv_t)\\
	\bJ^p_{\mcM}(\bv_t)\\
	\bJ^q_{\mcM}(\bv_t)
	\end{array}\right]~\textrm{and}~~
	\bJ_{\mcO}(\bv_t):=\left[\begin{array}{c}
	\bJ^p_\mcO(\bv_t)\\
	\bJ^q_\mcO(\bv_t)
	\end{array}\right].
	\end{align*}
\end{subequations}
Every $\bJ_{\mcM}(\bv_t)$ corresponds to $4M$ metering equations, and every $\bJ_{\mcO}(\bv_t)$ to $2O$ coupling equations. Having defined $\bJ_{\mcM}(\bv_t)$ and $\bJ_{\mcO}(\bv_t)$, the entire Jacobian matrix $\bJ\left(\{\bv_t\}\right)$ can be row-permuted as
\begin{equation}\label{eq:Jvv2}
\left[
\begin{array}{@{}c:c:c:c:c@{}} 
\mathbf{J}_{\mcM}(\bv_1) & \mathbf{0} &\mathbf{0} & \cdots & \mathbf{0}\\	
\mathbf{J}_{\mcO}(\bv_1) & -\mathbf{J}_{\mcO}(\bv_{2})  & \mathbf{0}  & \cdots & \mathbf{0} \\
\hdashline
\mathbf{0} & \mathbf{J}_{\mcM}(\bv_{2})  & \mathbf{0} & \cdots & \mathbf{0}\\
\mathbf{0} &\mathbf{J}_{\mcO}(\bv_{2}) & -\mathbf{J}_{\mcO}(\bv_{3})  &\cdots & \mathbf{0} \\
\hdashline
\mathbf{0} & \mathbf{0} & \mathbf{J}_{\mcM}(\bv_{3}) & \cdots & \mathbf{0}\\
\vdots & \vdots & \vdots & \ddots & \vdots\\
\mathbf{0} & \mathbf{0} & \mathbf{0} & \cdots & -\mathbf{J}_{\mcO}(\bv_T)\\
\hdashline
\mathbf{0} & \mathbf{0} & \mathbf{0} &  \cdots & \mathbf{J}_{\mcM}(\bv_T)
\end{array}
\right].
\end{equation}
This row-permuted version of $\bJ\left(\{\bv_t\}\right)$ will be denoted by $\tbJ\left(\{\bv_t\}\right)$, and has been obtained by interleaving block rows of metering and coupling equations. 

\begin{figure}[t]
\centering
\includegraphics[scale=0.29]{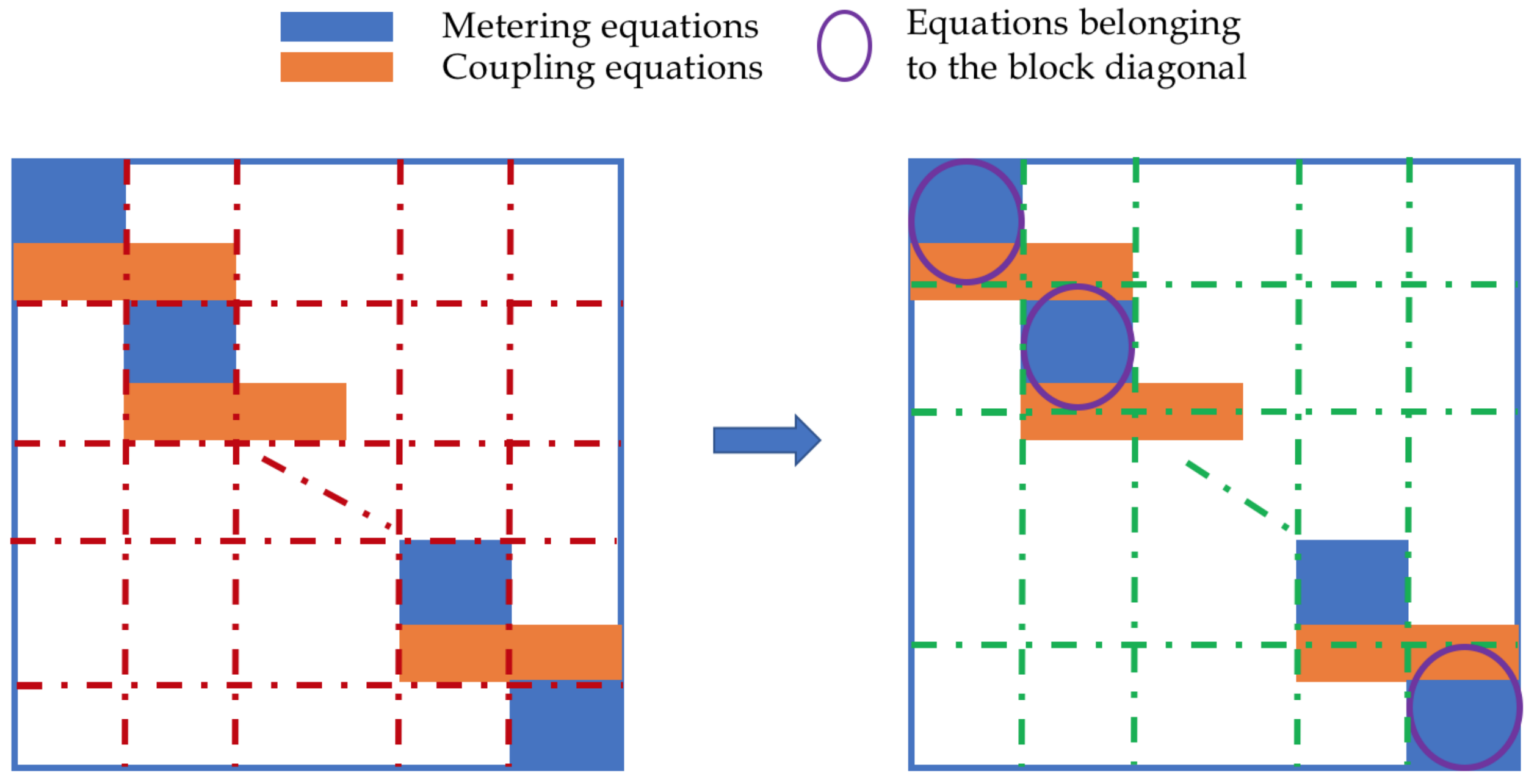}
\caption{\emph{Left:} Sparsity pattern of $\tbJ\left(\{\bv_t\}\right)$ [cf.~\eqref{eq:Jvv2}]. \emph{Right:} block tridiagonal $\tbJ\left(\{\bv_t\}\right)$ revealed after splitting each block row of coupling equations.}
\label{fig:block}
\end{figure}

Matrix $\tbJ\left(\{\bv_t\}\right)$ features the sparsity pattern of a block tridiagonal matrix. To reveal this structure, split each block row of coupling equations into two block rows. The top block row will be grouped with the previous block row of metering equations. The bottom block row will be grouped with the next block row of metering equations as in Fig.~\ref{fig:block}.

Focus now on the blocks lying on the main diagonal of $\tbJ\left(\{\bv_t\}\right)$. These blocks will be denoted by $\tbJ_t(\bv_t)$ for $t\in\mcT$. If for each $\tbJ_t(\bv_t)$, its columns can be perfectly matched to its rows, then a perfect bipartite matching for the entire $\tbJ\left(\{\bv_t\}\right)$ has been obtained. Then, Lemma~\ref{le:bipartite} guarantees that $\tbJ\left(\{\bv_t\}\right)$ and $\bJ\left(\{\bv_t\}\right)$ are generically full rank. 

Our goal is to assign coupling equations to blocks so that every block $\tbJ_t(\bv_t)$ enjoys a perfect bipartite matching. There are $2O(T-1)$ coupling equations to be assigned to $T$ blocks. A uniform allocation should assign $\frac{2O(T-1)}{T}$ coupling equations per block. With this allocation, block $t$ will have $4M$ metering equations and $\frac{2O(T-1)}{T}$ coupling equations over its $2(N+1)=2M+2O$ states in $\bv_t$. For a perfect bipartite matching to exist, we need $4M+\frac{2O(T-1)}{T}\geq 2M+2O$. 

The last requirement coincides with the necessary condition of \eqref{eq:cond2} for $T\geq 2$; but it is not enough: Every coupling equation can be assigned to exactly one between two specific blocks; see Fig.~\ref{fig:block}. For example, a coupling equation in the block row involving $\bJ_\mcO(\bv_2)$ and $-\bJ_\mcO(\bv_3)$ can be grouped either with $\bJ_\mcM(\bv_2)$ or $\bJ_\mcM(\bv_3)$. Partitioning the coupling equations into groups of $\frac{2O(T-1)}{T}$ while adhering to the latter requirement is the crux of the identifiability analysis. To allocate coupling equations, let us first define the bipartite grid graph $\mcG_b$.

\begin{definition}[Bipartite grid graph]\label{def:tG}
Consider the graph obtained from $\mcG$ upon maintaining only the edges between $\mcM$ and $\mcO$. Replicate the node set $\mcM$ to form $\mcM'$, and connect the nodes in $\mcM'$ to nodes in $\mcO$ by replicating the $\mcM$--$\mcO$ edges. The obtained bipartite graph will be denoted by $\mcG_b$. 
\end{definition} 

The identifiability of P2L relies on a matching in $\mcG_b$.

\begin{theorem}\label{th:main} 
If $\mcO$ can be partitioned into $\{\bmcO_k\}_{k=1}^{\ceil{T/2}}$ so that each one of them independently can be perfectly matched to $\mcM\cup\mcM'$ on $\mcG_b$, the Jacobian matrix $\bJ\left(\{\bv_t\}\right)$ related to the P2L task with phasor data is generically full rank.
\end{theorem}

In essence, Theorem~\ref{th:main} provides sufficient conditions for successful probing. The proof of Theorem~\ref{th:main} relies on two lemmas shown in the appendix: Lemma~\ref{le:block} provides sufficient conditions for the coupling equations assigned to block $t$, so that $\tbJ_t(\bv_t)$ enjoys a bipartite matching. Lemma~\ref{le:pattern} explains when these conditions can be met simultaneously for all $t\in\mcT$. The analysis uses the concept of a \emph{multi-set}. Different from a conventional set that contains unique elements, a multi-set is allowed to have multiple instances of elements. For example, we will override the definition of set union, so that $\{a,b\} \cup \{a,b\}$ does not yield $\{a,b\}$, but the multi-set $\{a,a,b,b\}$. 

\begin{lemma}\label{le:block}
Partition $\mcO$ into $\mcO_t\cup \bmcO_t$ so that $|\mcO\cup\mcO_t|=2O(T-1)/T$. Assume block $\tbJ_t(\bv_t)$ is assigned some coupling equations related to $\mcO\cup\mcO_t$. If the vertices in $\bmcO_t$ can be matched to the vertices in $\mcM\cup\mcM'$ on $\mcG_b$, the block $\tbJ_t(\bv_t)$ features a bipartite matching from its columns to its rows. 
\end{lemma}

\begin{lemma}\label{le:pattern}
Under the condition of Theorem~\ref{th:main}, the coupling equations for two successive blocks $\tbJ_t(\bv_t)$ with $t=2k-1$ and $t=2k$ share the same sparsity pattern of $\mcO\cup\mcO_k$ for $k=1,\ldots,\ceil{T/2}$.
\end{lemma}

\begin{figure}[t]
	\centering
	\includegraphics[scale=0.2]{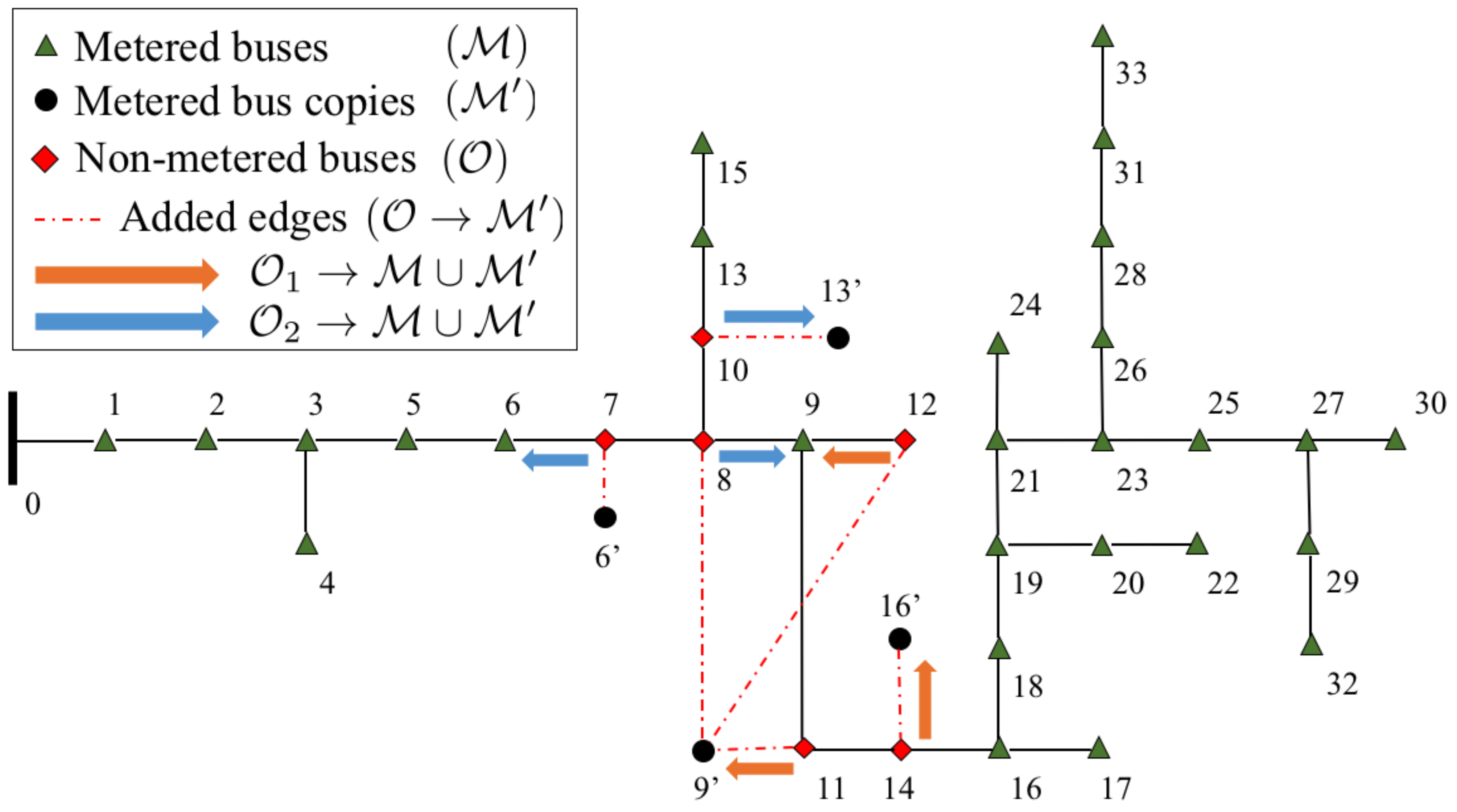}
	\caption{Matchings on the IEEE 34-bus grid for $T=4$ and $O=6$ for the P2L task with phasor data.}
	\label{fig:34busA}
\end{figure}

Theorem~\ref{th:main} follows as a direct consequence of Lemmas~\ref{le:block} and \ref{le:pattern}. To simplify the exposition, we will henceforth assume even $T$. To appreciate the conditions of Theorem~\ref{th:main}, examine the \emph{probing setup}, that is the placement of non-metered and probed buses, of Figure~\ref{fig:34busA}. The black circles denote the copies $\mcM'$ of nodes in $\mcM$, and the dashed red lines show the added edges from $\mcO$ to $\mcM'$. The operator needs to infer the loads at the $O=6$ non-metered buses marked by red diamonds. To study if probing this feeder over $T=4$ slots is successful, the set $\mcO$ has to be partitioned into two subsets $\mcO_1$ and $\mcO_2$, so that the buses of each subset are matched to buses in $\mcM\cup\mcM'$ on $\mcG_b$. The orange and blue arrows show precisely these matchings. If the feeder were to be probed over $T=2$ slots instead, probing would fail since buses $\{8,11,12\}$ cannot be uniquely matched to any buses in $\mcM\cup\mcM'$. 

\begin{algorithm}[t]
	\caption{Test for Successful Probing (phasor data)} \label{alg:phasor} 
	\begin{algorithmic}[1]
		\STATE Assign unit capacity to edges in $\mcG_b$ to define graph $\tmcG_b$.
		\STATE In $\tmcG_b$, add source node $n_s$, and connect it to all nodes in $\mcO$. These edges are assigned unit capacity.
		\STATE In $\tmcG_b$, add destination node $n_d$, and connect it to all nodes in $\mcM\cup\mcM'$.
		\STATE Initialize $T=2$.
		\WHILE{$T\leq T_{\max}$}
				\STATE The edges running between $\mcM\cup\mcM'$ and $n_d$ are assigned capacities of $T/2$.
				\STATE Run a max-flow problem between $n_s$ and $n_d$.
				\IF{obtained $n_s$--$n_d$ flow equals $O$} 
				\RETURN \emph{Probing setup is deemed successful for $T$.}
				\ELSE \STATE{$T:=T+2$}
				\ENDIF
				\RETURN \emph{Probing setup is deemed unsuccessful.}
		\ENDWHILE
	\end{algorithmic}
\end{algorithm}

As illustrated through this example, to check the condition of Theorem~\ref{th:main} for a particular $(\mcO,\mcM)$ probing setup, first one has to construct the bipartite graph $\mcG_b$ from $\mcG$. Then, given a number of probing actions $T$: \emph{i)} the set $\mcO$ has to be partitioned into the subsets $\{\bmcO_k\}_{k=1}^{T/2}$; and \emph{ii)} the nodes within each $\bmcO_k$ have to be mapped to the nodes in $\mcM\cup\mcM'$ on $\mcG_b$. Albeit these steps may seem computationally hard, they can be solved by a linear program as detailed in Algorithm~\ref{alg:phasor}. 

Given a probing setup, Algorithm~\ref{alg:phasor} finds the maximum flow between nodes $n_s$ and $n_d$ over graph $\tmcG_b$ constructed from $\mcG_b$. The edges in $\tmcG_b$ are organized in three layers: The edges of the first layer connect $n_s$ to $\mcO$ and have unit capacities. The edges of the second layer connect $\mcO$ to $\mcM\cup\mcM'$ and have unit capacities as well. The edges of the third layer connect $\mcM\cup\mcM'$ to $n_d$ and have capacities of $T/2$. This is to ensure that each node in $\mcM\cup\mcM'$ is mapped to at most $T/2$ nodes in $\mcO$ through the second layer. If the maximum $n_s$--$n_d$ flow equals $O$, all first-layer edges have been used to their capacity to map every node in $\mcO$ to exactly one node in $\mcM\cup\mcM'$. 

The max-flow problem can be solved using the Ford-Fulkerson algorithm, whose complexity scales linearly with the number of graph nodes and edges~\cite{Ford}. Moreover, if all edge capacities are integers, the algorithm finds an integral maximal flow. If the maximum $n_s$--$n_d$ flow is smaller than $O$, there is no matching for the tested $T$. Then, the edge capacities at the third layer can be increased and the process is repeated. Theorem~\ref{th:main} asserts that the chances of successful probing improve for larger $T$. This is because progressively smaller subsets of $\mcO$ need to be mapped to $\mcM\cup\mcM'$. Yet this gain in $T$ is limited by the bus placement $(\mcM,\mcO)$ as quantified next and shown in the appendix.

\begin{lemma}\label{le:2}
If $\delta_\mcM$ is the maximum node degree over $\mcM$ on $\tmcG_b$, a probing setup with phasor data cannot turn into successful beyond $T_{\max}=\delta_\mcM-1$. 
\end{lemma}

Lemma~\ref{le:2} implies that increasing $T$ beyond $T_{\max}$ has no hope in making probing successful for a specific placement, and Algorithm~\ref{alg:phasor} terminates with a negative answer. Once a $(\mcM,\mcO)$ placement is deemed successful, there are two questions to be answered: \emph{i)} how to select probing injections; and \emph{ii)} how to recover the non-metered loads. Both questions along with numerical tests are deferred to Part II.

\section{Identifiability of P2L with Non-phasor Data}\label{sec:non-phasor}
Since PMUs have limited penetration in distribution grids, requiring voltage phasor data at probing buses may be unrealistic. This section studies probing with non-phasor data.

\begin{figure}[t]
	\centering
	\includegraphics[scale=0.2]{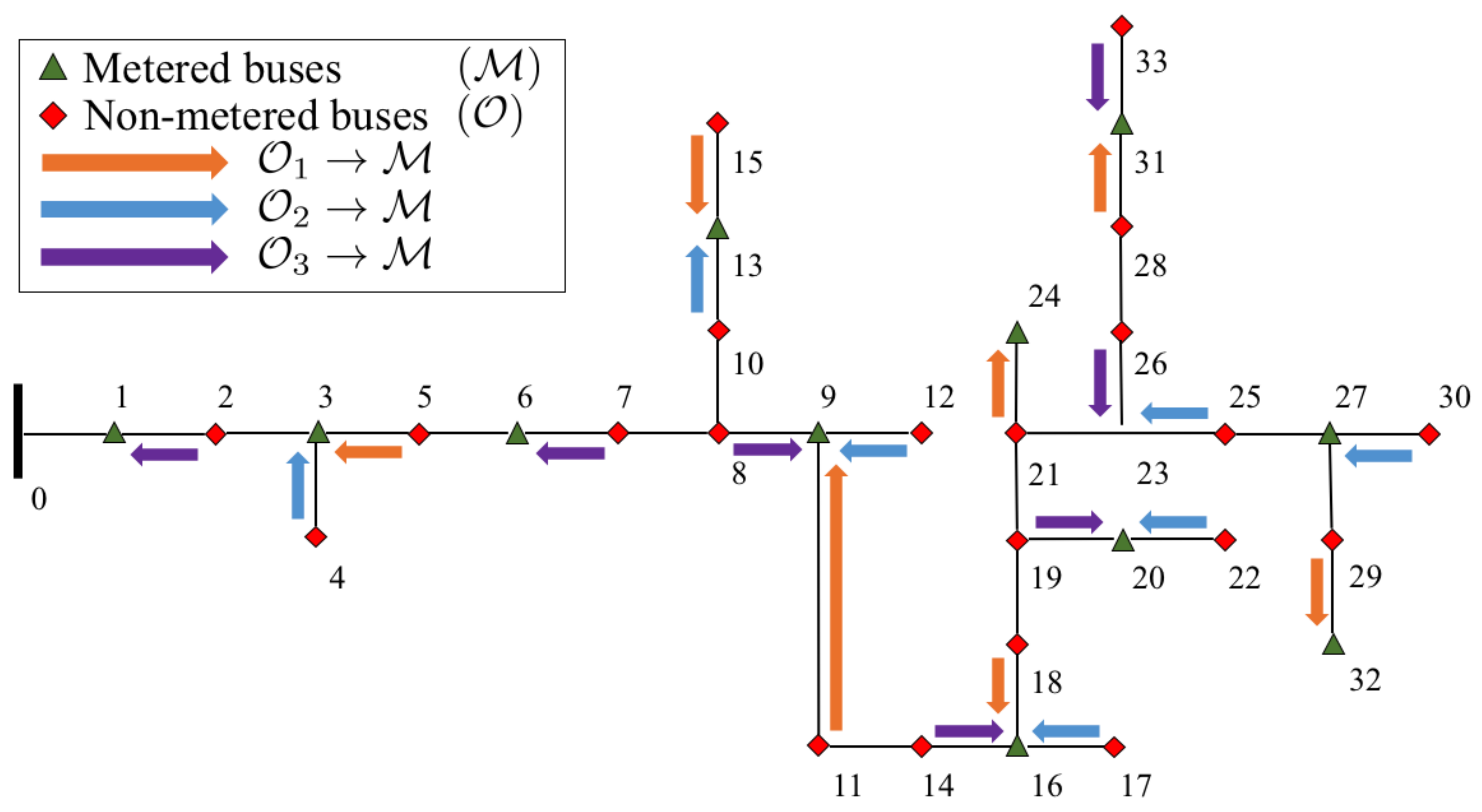}
	\caption{Matchings on the IEEE 34-bus grid for $T=6$ and $O=21$ for the P2L task with non-phasor data.}
	\label{fig:34busB}
\end{figure}

\begin{definition}[P2L task with non-phasor data]\label{def:P2L2}
Given $\mathbf{Y}$ and probing data $(u_n^t, p_n^t,q_n^t)$ for $n\in\mcM$ and $t\in\mcT$, the P2L task entails solving the equations in \eqref{eq:PF1u} and \eqref{eq:PF1pq} for $t\in\mcT$, jointly with the coupling equations in \eqref{eq:PFC}. 
\end{definition}

A simple count of equations and unknowns dictates $M\geq\frac{2O}{T}$,  which is clearly more restrictive than \eqref{eq:cond2}. We next provide a sufficient condition under which this task is solvable.

\begin{theorem}\label{th:main2}
If $\mcO$ can be partitioned into $\{\bmcO_k\}_{k=1}^{\ceil{T/2}}$ such that each one of them independently can be perfectly matched to $\mcM$ on $\mcG$, the Jacobian matrix $\bJ\left(\{\bv_t\}\right)$ related to the P2L task with non-phasor data is generically full rank.
\end{theorem}

Dropping the voltage angle metering equations, matrix ${\bJ}_{\mcM}(\bv_t)$ in \eqref{eq:Jvv2} is replaced by 
	\begin{subequations}
		\begin{align*}
			{\bJ}_{\mcM}(\bv_t):= \left[\begin{array}{c}
				\bJ^u_\mcM(\bv_t)\\
				\bJ^p_{\mcM}(\bv_t)\\
				\bJ^q_{\mcM}(\bv_t)
			\end{array}\right].
		\end{align*}
	\end{subequations}
Similar to Theorem~\ref{th:main}, it is not hard to see that the nodes in $\bmcO_t$ have to be matched to the nodes in $\mcM$, rather than $\mcM \cup \mcM'$.

Consider for example the probing setup of Figure~\ref{fig:34busB}. To infer the loads at $O=21$ non-metered buses with $T=6$ probing slots, the set $\mcO$ has to be partitioned into three subsets $\mcO_1$, $\mcO_2$, and $\mcO_3$, so that the buses of each subset are matched to $\mcM$. The orange, blue, and purple arrows in the figure show these matchings. Because non-metered buses are divided into three subsets, up to three non-metered buses can be matched to the same probed bus. For example, buses $\{14,17,18\}$ are all matched to the probed bus $16$. Probing the same feeder over $T=2$ or $T=4$ rather than $T=6$ slots would fail.

\begin{algorithm}[t]
	\caption{Test for Successful Probing (non-phasor data)} \label{alg:non-phasor} 
	\begin{algorithmic}[1]
		\STATE Connect $n_s$ to all nodes in $\mcO$ with unit-capacity edges.
		\STATE Connect $\mcO$ to $\mcM$ based on $\mcG$ with unit-capacity edges.
		\STATE Connect all nodes in $\mcM$ to $n_d$.
		\STATE Initialize $T=2$.
		\WHILE{$T\leq T_{\max}$}
				\STATE Assign capacity $T/2$ to edges between $\mcM$ and $n_d$.
				\STATE Run a max-flow problem between $n_s$ and $n_d$.
				\IF{obtained $n_s$--$n_d$ flow equals $O$} 
				\RETURN \emph{Probing setup is deemed successful for $T$.}
				\ELSE \STATE{$T:=T+2$}
				\ENDIF
				\RETURN \emph{Probing setup is deemed unsuccessful.}
		\ENDWHILE
	\end{algorithmic}
\end{algorithm}

The condition of Theorem~\ref{th:main2} can be easily tested by Algorithm~\ref{alg:non-phasor} and up to the value of $T_{\max}$ provided next. 

\begin{corollary}\label{co:1}
If $\delta_\mcM$ is the maximum degree of the nodes in $\mcM$ on the graph constructed by Alg.~\ref{alg:non-phasor}, a probing setup with non-phasor data cannot turn into successful beyond $T_{\max}=2(\delta_\mcM-1)$.
\end{corollary}  

Corollary~\ref{co:1} is proved as part of the proof of Lemma~\ref{le:2}. Compared to Theorem~\ref{th:main2}, the condition of Theorem~\ref{th:main} provided more flexibility towards attaining a bipartite matching since probed buses can be used twice. If a probing setup is successful for non-phasor data, it is also successful for phasor data. Interestingly, the matchings in Theorems~\ref{th:main} and \ref{th:main2} depend solely on the sparsity pattern of $\bG$ and the probing setup, so the claims here apply to even meshed (e.g., multiphase) grids.



\begin{figure}[t]
	\centering
	\includegraphics[scale=0.34]{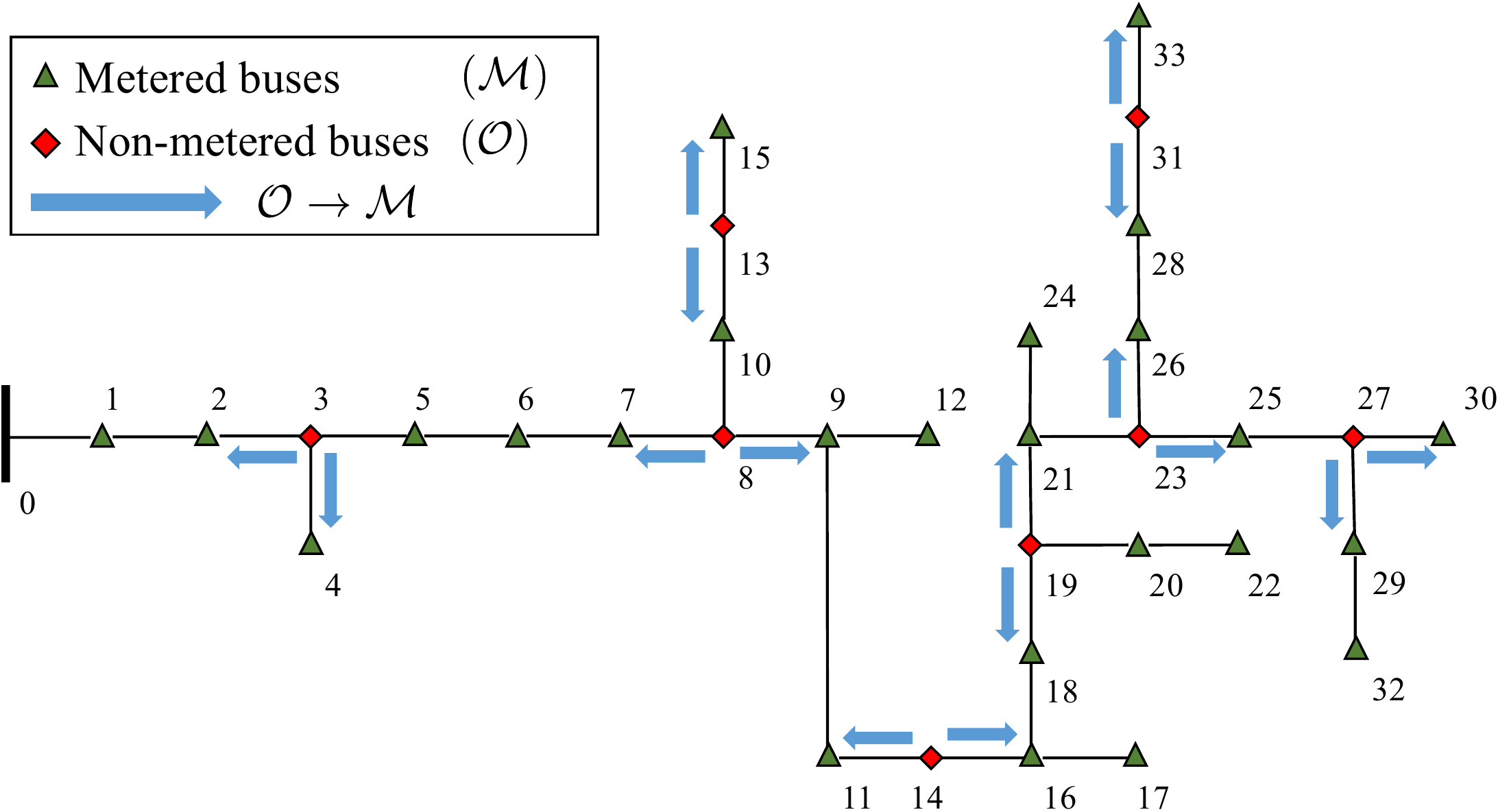}
	\caption{Matchings on the IEEE 34-bus grid for single-slot probing $(T=1)$ with non-phasor data and $O=8$. Single-slot probing waives assumption \emph{a2)} on exclusively constant-power loads.}
	\label{fig:34busD}
\end{figure}

\section{Single-Slot Probing} \label{sec:SSPF}
The analysis so far depends on assumption \emph{a2}) of constant-power loads. Under the ZIP load model of \eqref{eq:ZIP0}, the coupling equations in \eqref{eq:PFC} are no longer valid, and thus, the metering equations decouple across $\mcT$. Can the non-metered loads $p_n+jq_n$ for $n\in\mcO$ still be recovered upon collecting data on $\mcM$? This answer can be on the affirmative with \emph{single-slot probing}, that is $T=1$. Leveraging the tools of Sections~\ref{sec:phasor} and \ref{sec:non-phasor}, we next study the observability of single-slot probing. The ensuing two results (proven in the appendix) provide conditions for successful load recovery using (non)-phasor data.




\begin{theorem}\label{th:main3}
If each bus in $\mcO$ can be matched to one unique bus in $\mcM$ on $\mcG$, the Jacobian matrix $\bJ(\bv_1)$ related to single-slot probing $(T=1)$ with phasor data has full generic rank.
\end{theorem}

\begin{theorem}\label{th:main4} 
If each bus in $\mcO$ can be matched to two unique buses in $\mcM$ on $\mcG$, the Jacobian $\bJ(\bv_1)$ related to single-slot probing $(T=1)$ with non-phasor data has full generic rank.
\end{theorem}

The conditions of Th.~\ref{th:main3} and \ref{th:main4} can be tested by Algorithm~\ref{alg:non-phasor} by fixing $T=2$ and $T=1$, respectively. Figure~\ref{fig:34busD} shows a successful placement per Theorem~\ref{th:main4}. 

If the conditions of Th.~\ref{th:main3} and \ref{th:main4} are met, the non-metered loads $p_n+jq_n$ for $n\in\mcO$ can be recovered using single-slot probing, regardless if these loads are constant-power or not. However, the operator may also want to estimate their ZIP parameters in \eqref{eq:ZIP0}. Estimating these parameters directly with multi-slot probing becomes complicated. Instead, one could adopt multi-slot probing in a two-step process as follows: First, the feeder is probed over $\mcT$ with $|\mcT|=T>3$. Under Th.~\ref{th:main3} and \ref{th:main4}, the operator obtains estimates $(\hat{u}_n^t,\hat{p}_n^t,\hat{q}_n^t)$ for all non-metered buses $n\in\mcO$ and $t\in\mcT$. Secondly, the ZIP parameters for active load $n$ can be estimated through the least-squares (LS) fit
\begin{equation} \label{eq:ZIP}
[\hat{\alpha}_{p_n}~\hat{\beta}_{p_n}~\hat{\gamma}_{p_n}]^\top
:=(\hbU_n^{\top}\hbU_n)^{-1}\hbU_n^{\top}\hbp_n
\end{equation}
where $\hbp_n:=[\hat{p}_n^1~\ldots~\hat{p}_n^T]^\top$ and the $t$-th row of matrix $\hbU_n$ is $[(\hat{u}_n^t)^2~\hat{u}_n^t~1]$ for $t=1,\ldots,T$. Similar LS fits can be performed for the reactive ZIP load parameters. A major concern here is that all entries of the Vandermonde matrix $\hbU_n$ are close to unity in compliance with voltage regulation. For $T=3$, the determinant of $\hbU_n$ is calculated as $(\hat{u}_n^1-\hat{u}_n^2)(\hat{u}_n^1-\hat{u}_n^3)(\hat{u}_n^2-\hat{u}_n^3)$~\cite{HornJohnson}, which yields $|\hbU_n|=-2\cdot 10^{-3}$ even for bus voltages as widely spread as $\hat{u}_n^1=0.9$, $\hat{u}_n^2=1.0$, and $\hat{u}_n^3=1.1$. This reveals that the task of estimating ZIP parameters from voltage/power data is \emph{ill-posed}. This is germane to the task itself rather than the method (here probing) used to collect the data. 

Finally, note that Th.~\ref{th:main3} and \ref{th:main4} hold even when data are not collected via probing, e.g., smart meter data. Therefore, our observability analysis covers the general setup where voltage and active/reactive power data or specifications are given only for $\mcM$. Similar conditions were derived in \cite{BKVGS17}, but were confined to radial grids. 


\section{Conclusions}\label{sec:conclusions}
The novel technique of intentionally probing an electric grid using inverters to recover non-metered loads has been put forth in the first part of this two-part work. The technique leverages the actuation capabilities of smart inverters, the data collected at probed buses, and the stationarity of non-metered loads, to formulate a power flow problem coupled over multiple times. Sufficient conditions that can be easily verified by solving a max-flow problem on a grid graph have been provided to test if a probing placement is successful. Beyond probing, the pertinent task of finding loads using data from a subset of buses has also been cast as a special case. Assuming a probing setup satisfies these conditions, Part II explains how inverter probing setpoints can be designed to improve load estimation accuracy, and provides numerical tests on the IEEE 34-bus feeder using a semidefinite program relaxation.

\appendix 
\begin{IEEEproof}[Proof of Lemma~\ref{le:block}]
It can be easily verified that the sparsity patterns of $\mathbf{J}^u(\bv_t)$ and $\mathbf{J}^{\theta}(\bv_t)$ coincide with the sparsity pattern of $[\mathbf{I}_{N+1}~\mathbf{I}_{N+1}]$. The sparsity patterns of $\bJ^p(\bv_t)$ and $\bJ^{q}(\bv_t)$ coincide with the sparsity pattern of $[\mathbf{G}~\mathbf{G}]$ where $\bG$ is the bus conductance matrix; see \cite[Table 3.2]{ExpConCanBook}. From \eqref{eq:Jvv2}, the sparsity pattern of $\tbJ_{t}(\bv_t)$ is 
\begin{equation} \label{eq:split1}
\begin{bmatrix}
\bI_{\mcM,\mcN^+} & \bI_{\mcM,\mcN^+}\\
\bI_{\mcM,\mcN^+} & \bI_{\mcM,\mcN^+}\\
\bG_{\mcM,\mcN^+} & \bG_{\mcM,\mcN^+}\\
\bG_{\mcM,\mcN^+} & \bG_{\mcM,\mcN^+}\\
\bG_{\mcO,\mcN^+} & \bG_{\mcO,\mcN^+}\\
\bG_{\mcO_t,\mcN^+} & \bG_{\mcO_t,\mcN^+}
\end{bmatrix}
\end{equation}
where the first block row relates to voltage magnitudes; the second to voltage angles; the third and fourth to probing injections; while the fifth and sixth to coupled injections. 

To create a bipartite matching for block $\tbJ_{t}(\bv_t)$, unfold the sparsity pattern in \eqref{eq:split1} column-wise using $\mcN^+=\mcM\cup \mcO$ as
\begin{equation} \label{eq:split2}
\begin{bmatrix}
\textcolor{blue}{\fbox{$\bI_{\mcM,\mcM}$}} & \bI_{\mcM,\mcO} & \bI_{\mcM,\mcM} & \bI_{\mcM,\mcO} \\
\bI_{\mcM,\mcM} & \bI_{\mcM,\mcO} & \textcolor{blue}{\fbox{$\bI_{\mcM,\mcM}$}} & \bI_{\mcM,\mcO} \\
\bG_{\mcM,\mcM} & \bG_{\mcM,\mcO} & \bG_{\mcM,\mcM} & \textcolor{red}{\dbox{$\bG_{\mcM,\mcO}$}}\\
\bG_{\mcM,\mcM} & \bG_{\mcM,\mcO} & \bG_{\mcM,\mcM} & \textcolor{red}{\dbox{$\bG_{\mcM,\mcO}$}}\\
\bG_{\mcO,\mcM} & \textcolor{blue}{\fbox{$\bG_{\mcO,\mcO}$}} & \bG_{\mcO,\mcM} & \bG_{\mcO,\mcO}\\
\bG_{\mcO_t,\mcM} & \bG_{\mcO_t,\mcO} & \bG_{\mcO_t,\mcM} & \textcolor{red}{\dbox{$\bG_{\mcO_t,\mcO}$}}
\end{bmatrix}.
\end{equation}
The first block column in \eqref{eq:split2} relates to variables $\{v_{r,n}^t\}_{n\in\mcM}$, and can be matched to the first block row via $\bI_{\mcM,\mcM}$. Similarly, the third block column relates to variables $\{v_{i,n}^t\}_{n\in\mcM}$, and can be matched to the second block row. The second block column relates to variables $\{v_{r,n}^t\}_{n\in\mcO}$, and can be matched to the fifth block row via the main diagonal of $\bG_{\mcO,\mcO}$. 

To achieve a bipartite matching, the fourth block column related to variables $\{v_{i,n}^t\}_{n \in \mcO}$ has to be matched to the union of the third, fourth, and sixth block rows. Lacking a simple diagonal matching now, we leverage the sparsity pattern of $\bG$. It suffices to match the column nodes in $\mcO$ to the row nodes in $\mcM \cup \mcM \cup \mcO_t$. Because $\mcO=\mcO_t\cup \bmcO_t$, the column nodes $\mcO_t$ can be matched to the row nodes $\mcO_t$ via some diagonal entries of $\bG_{\mcO_t,\mcO}$. Then, the column nodes $\bmcO_t$ have to be matched to the row nodes $\mcM \cup \mcM$. This can be accomplished based on the hypothesis of this Lemma, thus completing its proof.
\end{IEEEproof}

\begin{IEEEproof}[Proof of Lemma~\ref{le:pattern}]
The pair of blocks $\tbJ_{2k-1}(\bv_{2k-1})$ and $\tbJ_{2k}(\bv_{2k})$ will be jointly indexed by $k$. Define also
\begin{equation}\label{eq:Rk}
\mcR_k:=\bigcup_{\tau=1}^k \bmcO_{\tau}.
\end{equation}
In addition to the claim of this lemma, we will also prove that when passing from pair $k-1$ to pair $k$, a set of coupling equations represented by $\mcR_{k-1}\cup\mcR_{k-1}$ have not been assigned to block $2k-2$, and are free to be assigned to block $2k-1$.

Proving by induction, we start with the base case. The pair indexed by $k=1$ consists of $\tbJ_1(\bv_1)$ and $\tbJ_2(\bv_2)$. The active and reactive equations coupling these two blocks can be represented by $\mcO\cup\mcO=\mcO_1\cup\mcO_1\cup\bmcO_1\cup\bmcO_1$. Let us assign $\mcO_1\cup\mcO_1\cup\bmcO_1$ to block 1. With this assignment, the coupling equations for $\tbJ_1(\bv_1)$ get the sparsity pattern of $\mcO\cup\mcO_1$. The remaining coupling equations in $\bmcO_1$ are assigned to block 2. 

Block 2 shares with block 3 the coupling equations $\mcO\cup\mcO$, which are again expressed as $\mcO_1\cup\mcO_1\cup\bmcO_1\cup\bmcO_1$. From this new set of coupling equations, assign $\mcO_1\cup\mcO_1$ to block 2. Hence, the coupling equations for $\tbJ_2(\bv_2)$ have the sparsity pattern of $\mcO_1\cup\mcO_1\cup\bmcO_1=\mcO\cup\mcO_1$. The unused coupling equations are represented by $\bar{\mcO}_1\cup\bar{\mcO}_1=\mcR_1\cup\mcR_1$.

Suppose the claim holds for the block pair $k-1$. It is next shown that the claim holds for the block pair $k$ consisting of blocks $2k-1$ and $2k$. Starting with the odd block $2k-1$, the unused equations $\mcR_{k-1}\cup \mcR_{k-1}$ that couple blocks $2k-2$ and $2k-1$ are assigned to block $2k-1$. Block $2k-1$ is also coupled to block $2k$ via $\mcO\cup\mcO$ equations, which can be expressed as $\mcO_k\cup\mcO_k\cup\bmcO_k\cup\bmcO_k$. The key point here is that by the definition of $\mcR_{k-1}$ and because $\bmcO_k$'s are mutually exclusive by the hypothesis of this lemma, it holds that 
\[\mcR_{k-1} \cap \bmcO_{k} = \emptyset~\text{and}~\mcR_{k-1} \subset \mcO,~\text{so that}~\mcR_{k-1} \subseteq \mcO_{k}.\]
Therefore, the set $\mcO_k$ can be partitioned into $\mcR_{k-1}$ and $\mcO_{k} \setminus \mcR_{k-1}$. From the equations coupling blocks $2k-1$ and $2k$, the equations $(\mcO_{k}\setminus \mcR_{k-1})\cup(\mcO_{k}\setminus \mcR_{k-1}) \cup \bmcO_{k}$ are assigned to block $2k-1$. In this way, the coupling equations for block $2k-1$ have the sparsity pattern
\[\underbrace{\mcR_{k-1}\cup \mcR_{k-1}}_{\text{with block $2k-2$}}\cup \underbrace{(\mcO_{k}\setminus \mcR_{k-1})\cup(\mcO_{k}\setminus \mcR_{k-1}) \cup \bmcO_{k}}_{\text{with block $2k$}}=\mcO\cup\mcO_k.\]
The unused equations coupling blocks $2k-1$ and $2k$ are $\mcR_{k-1}\cup\mcR_{k-1}\cup\bmcO_k$.

Moving to block $2k$ of pair $k$, the unused equations $\mcR_{k-1}\cup\mcR_{k-1}\cup\bmcO_k$ coupling block $2k$ with block $2k-1$ are assigned to block $2k$. Block $2k$ is also coupled with block $2k+1$ through $\mcO\cup\mcO=\mcO_k\cup\mcO_k\cup\bmcO_k\cup\bmcO_k$. From this new set of coupling equations, assign equations $(\mcO_{k}\setminus \mcR_{k-1})\cup(\mcO_{k}\setminus \mcR_{k-1})$ to block $2k$. Hence, the coupling equations assigned to block $2k-1$ have the sparsity pattern
\[\underbrace{\mcR_{k-1}\cup \mcR_{k-1}\cup\bmcO_k}_{\text{with block $2k-1$}}\cup \underbrace{(\mcO_{k}\setminus \mcR_{k-1})\cup(\mcO_{k}\setminus \mcR_{k-1})}_{\text{with block $2k+1$}}=\mcO\cup\mcO_k.\]
The unused equations coupling blocks $2k$ and $2k+1$ are
\[\mcR_{k-1}\cup\mcR_{k-1}\cup\bmcO_k\cup\bmcO_k=\mcR_k\cup\mcR_k.\]

For the last block pair, the coupling equations already assigned to block $T-1$ have the sparsity pattern $\mcO\cup\mcO_{T/2}$. The last block $T$ differs from the previous blocks as it only gets the $\mcR_{T/2-1}\cup\mcR_{T/2-1}\cup\bar{\mcO}_{T/2}$ unused coupling equations between blocks $T-1$ and $T$. Under the condition of Th.~\ref{th:main}, $\mcR_{T/2-1}\cup\bar{\mcO}_{T/2}=\mcO$ and because $\mcO={\mcO}_{T/2}\cup\bar{\mcO}_{T/2}$, we also have $\mcR_{T/2-1}=\mcO_{T/2}$. Hence, the sparsity pattern of the last block is also given by $\mcO\cup\mcO_{T/2}$. Since every pair of successive blocks has the same structure, the diagonal blocks of $\tbJ\left(\{\bv_t\}\right)$ will exhibit $\ceil{T/2}$ distinct sparsity patterns.
\end{IEEEproof}

\begin{IEEEproof}[Proof of Lemma~\ref{le:2}]
Consider node $m\in\mcM$ in $\tmcG_b$ with degree $\delta_m$. This node is connected to node $n_d$ via an edge having capacity $T/2$, and to $\delta_m-1$ nodes in $\mcO$ via unit-capacity edges. The maximum flow that can pass through the second-layer edges to $m$ is $\delta_m-1$. This flow will be funneled through edge $(m,n_d)$. Then, there is no advantage for this edge to have capacity larger than $\delta_m-1$, so that $T/2\leq \delta_m-1$. Considering all $m\in\mcM\cup\mcM'$, there is no point in testing for values of $T$ beyond $T\leq 2(\delta_\mcM-1)$.

The bound can be improved, since the previous argument assumed that all $\delta_m-1$ edges between $\mcO$ and $m\in\mcM$ have reached their capacity. That will not happen since the $\mcO$ nodes adjacent to $m$ on the feeder, can be shared between $m$ and its copy $m'\in\mcM'$ on $\tmcG_b$. Hence, the flow passing jointly through $m$ and $m'$ cannot exceed $\delta_m-1$. Then, the capacity of edge $(m,n_d)$ plus the capacity of edge $(m',n_d)$ can be safely limited to $\delta_m-1$, implying $T\leq \delta_m-1$. 
\end{IEEEproof}

\begin{IEEEproof}[Proof of Theorem~\ref{th:main3}]
The sparsity pattern of $\bJ(\bv_t)$ can be derived from \eqref{eq:split1}--\eqref{eq:split2} by eliminating the blocks related to coupling equations 
	\begin{equation} \label{eq:split3}
	\begin{bmatrix}
	\textcolor{blue}{\fbox{$\bI_{\mcM,\mcM}$}} & \bI_{\mcM,\mcO} & \bI_{\mcM,\mcM} & \bI_{\mcM,\mcO} \\
	\bI_{\mcM,\mcM} & \bI_{\mcM,\mcO} & \textcolor{blue}{\fbox{$\bI_{\mcM,\mcM}$}} & \bI_{\mcM,\mcO} \\
	\bG_{\mcM,\mcM} & \textcolor{red}{\dbox{$\bG_{\mcM,\mcO}$}} & \bG_{\mcM,\mcM} & \bG_{\mcM,\mcO}\\
	\bG_{\mcM,\mcM} & \bG_{\mcM,\mcO} & \bG_{\mcM,\mcM} & \textcolor{red}{\dbox{$\bG_{\mcM,\mcO}$}}
	\end{bmatrix}.
	\end{equation}
Variables  $\{v_{r,n}^t\}_{n\in\mcM}$ corresponding to the first block column and $\{v_{i,n}^t\}_{n\in\mcM}$ to the third block column can be matched respectively to the first and second block row via $\mathbf{I}_{\mcM,\mcM}$. To complete the bipartite matching, the second and fourth block columns (variables  $\{v_{r,n}^t\}_{n\in\mcO}$ and  $\{v_{i,n}^t\}_{n\in\mcO}$) can be matched to the third and fourth block rows, accordingly. Hence, $\bJ(\bv_t)$ is generically full rank if there exists a perfect matching in $\bG_{\mcM,\mcO}$, that is every non-metered node in $\mcO$ is mapped to a unique node in $\mcM$; see also Lemma~\ref{le:bipartite}. 
\end{IEEEproof}

\begin{IEEEproof}[Proof of Theorem~\ref{th:main4}]	
Given non-phasor data, the second block row related to voltage angles in \eqref{eq:split3} is dropped. Following the arguments to the proof for Theorem~\ref{th:main3}, matrix $\bJ(\bv_t)$ can be shown to be generically full rank if there exists a perfect matching in $\bG_{\mcM,(\mcO \cup \mcO)}$, that is every node in $\mcO$ is mapped to two unique nodes in $\mcM$; see also Lemma~\ref{le:bipartite}.  
\end{IEEEproof}

\balance
\bibliographystyle{IEEEtran}
\bibliography{myabrv,power}

\end{document}